\documentclass[12pt, A4paper]{amsart}

  \usepackage{latexsym,enumerate,xspace}
  \usepackage[all]{xy}
  \usepackage{amsthm}
 \numberwithin{equation}{section}
  \usepackage{color}
  \usepackage{hyperref}

\newcommand{\K}{\ensuremath{\cK}\xspace}
\renewcommand{\L}{\ensuremath{\mathcal L}\xspace}
\renewcommand{\epsilon}{\varepsilon}
\renewcommand{\phi}{\varphi}

\def\two{\ensuremath{\mathbf{2}\xspace}}
\def\op{^{\textnormal{op}}}

\def\cK{\mathcal K}

\def\1c#1{\stackrel{#1}{\to}}

\def\Mnd{\mathsf{Mnd}}
\def\Cmd{\mathsf{Cmd}}
\def\WDL{\mathsf{Wdl}}
\def\MDL{\mathsf{Mdl}}
\def\x{\times}
\def\Bim{\mathsf{Bim}}
\def\Cat{\ensuremath{\mathsf{Cat}}\xspace}

\def\Alg{\mathsf{alg}}
\def\Calg{\mathsf{coalg}}

\def\mone{\pi}
\def\mtwo{\sigma}
\def\monesec{\check{\pi}}
\def\mtworetr{\hat{\sigma}}
\def\mthree{\theta}

  \newtheorem{proposition}{Proposition}[section]
  \newtheorem{lemma}[proposition]{Lemma}

  \newtheorem{theorem}[proposition]{Theorem}

  \theoremstyle{definition}

  \theoremstyle{remark}
  \newtheorem{remark}[proposition]{Remark}

  \newcounter{c}
  
  \newcommand{\etyk}[1]{\vspace{-7.4mm}$$\begin{equation}\Label{#1}
  \addtocounter{c}{1}}
  \renewcommand{\]}{\ifnum \value{c}=1 $$\else \end{equation}\fi}
\begin{document}

 \title{On the 2-categories of weak distributive laws}

 \author{Gabriella B\"ohm}
 \address{Research Institute for Particle and Nuclear Physics, Budapest,
 \newline\indent H-1525
 Budapest 114, P.O.B.\ 49, Hungary}
 \email{G.Bohm@rmki.kfki.hu}

 \author{Stephen Lack}
 \address{School of Computing and Mathematics University of Western Sydney,
 \newline\indent 
 Locked Bag 1797 Penrith South DC NSW 1797, Australia and
 \newline \indent 
 Mathematics Department Macquarie University, NSW 2109 Australia.}
 \email{s.lack@uws.edu.au; steve.lack@mq.edu.au}

 \author{Ross Street}
 \address{Mathematics Department Macquarie University, NSW 2109 Australia.}
 \email{ross.street@mq.edu.au}
\dedicatory{Dedicated to Mia Cohen on the occasion of her retirement} 
 \date{Sep 2010}
 \begin{abstract}
A weak mixed distributive law (also called weak entwining structure
\cite{CaeDeG}) in a 2-category consists of a monad and a comonad, together with
a 2-cell relating them in a way which generalizes a mixed distributive law due
to Beck.     
We show that a weak mixed distributive law can be described as a compatible
pair of a monad and a comonad, in 2-categories extending, respectively, the
2-category of comonads and the 2-category of monads in \cite{Street}. 
Based on this observation, we define a 2-category whose 0-cells are weak
mixed distributive laws. 
In a 2-category $\cK$ which admits Eilen\-berg-Moore constructions both for
monads and comonads, and in which idempotent 2-cells split, we construct a
fully faithful 2-functor from this 2-category of weak mixed distributive laws
to $\cK^{2\x 2}$.  
\end{abstract}
  \maketitle


\section*{Introduction}

Distributive laws -- between two monads; between two comonads; or between a
monad and a comonad in any bicategory  (the latter known as the `mixed' case)
-- were discussed by Beck in \cite{Beck}. In Hopf algebra theory  the mixed
case was  introduced in \cite{BrzMaj} by Brzezi\'nski and Majid in the
particular bicategory $\Bim$ of Algebras; Bimodules; Bimodule Maps, under the
name `entwining structure', as a tool for unifying various Hopf type modules.   

For any 2-category \K, there is a 2-category $\Mnd(\K)$ of monads in \K for
which a monad in $\Mnd(\K)$ is the same thing as two monads in \K with a
distributive law between them \cite{Street}. Similarly, a comonad in
$\Mnd(\K)$ is the same thing as a mixed distributive law. Dually, there is a
2-category $\Cmd(\K)=\Mnd(\K_*)_*$, where $(-)_*$ denotes the vertical
opposite of a  2-category (the superscript ``co'' is also often used). A monad
in $\Cmd(\K)$ is once again a mixed distributive law, while a comonad in
$\Cmd(\K)$ is the same as two comonads with a distributive law between them.  

We identify the isomorphic 2-categories $\Mnd(\Cmd(\K))$ and $\Cmd(\Mnd(\K))$,
and write each as $\MDL(\K)$, the 2-category of mixed distributive laws. A
typical object will be written as  $(K,t,c,\lambda)$, where $K$ is the
underlying object, $t$ the monad, $c$ the comonad, and $\lambda:tc\to ct$ the
2-cell between them giving the distributive law. In the case where $t$ has a
right adjoint $d$, the monad structure on $t$ induces a comonad structure on
$d$, and mixed distributive laws $tc\to ct$ are in bijection with distributive
laws $cd\to dc$. 

Distributive laws play a key role in the description of liftings of monads and
comonads \cite{Street}, \cite{PowWat}. If the 2-category $\cK$ admits
Eilen\-berg-Moore constructions for monads; that is, the fully faithful
inclusion 2-functor $I:\cK \to \Mnd(\cK)$ possesses a right 2-adjoint $\Alg$,
then there is a fully faithful 2-functor from $\Mnd(\K)$ to the arrow
2-category $\K^\two$ (i.e. the 2-category of 2-functors from the interval
2-category $0\to 1$ to $\cK$). It sends a monad $(K,t)$ to the forgetful map
$u^t:K^t\to K$, seen as an object of $\K^\two$. The 2-functor $\Cmd(\Alg)$
takes the comonad $((K,t),(c,\lambda))$ in $\Mnd(\cK)$ to a comonad
$(K^t=\Alg(K,t),{\overline c}:=\Alg(c,\lambda))$ in $\cK$, which is a lifting
of the comonad $(K,c)$ to the Eilen\-berg-Moore object $K^t$ of
$(K,t)$.    

Similarly, if \K admits Eilen\-berg-Moore constructions for comonads; that is,
the fully faithful inclusion 2-functor $I_*:\K\to\Cmd(\K)$ has a right
2-adjoint $\Calg$, then there is a fully faithful 2-functor
$\Cmd(\K)\to\K^\two$ sending a comonad $(K,c)$ to the forgetful map
$u^c:K^c\to K$. Once again, the 2-functor $\Mnd(\Calg)$ takes the monad
$((K,c),(t,\lambda))$ in $\Cmd(\K)$ to a monad
$(K^c=\Calg(K,c),\overline{t}=\Calg(t,\lambda))$ which is a lifting of the
monad $(K,t)$ to the Eilen\-berg-Moore object $K^c$ of $(K,c)$.   

If \K admits Eilen\-berg-Moore constructions for both monads and comonads then
there is a commutative diagram of fully faithful 2-functors
$$
\xymatrix{
\K \ar[r]^I \ar[d]_{I_*} & \Mnd(\K) \ar[d]^{\Mnd(I_*)} \\
\Cmd(\K) \ar[r]^-{\Cmd(I)} & \MDL(\K) }
$$
and since both $\Alg\circ\Mnd(\Calg)$ and $\Calg\circ\Cmd(\Alg)$ are right
adjoint to the common diagonal in this last displayed diagram, they are
naturally isomorphic, sending an object $(K,t,c,\lambda)$ of $\MDL(\K)$ to
$(K^c)^{\overline{t}}$, respectively to $(K^t)^{\overline{c}}$. We shall
sometimes write $K^{(t,c)}$ for this common value, although really $\lambda$
should be included in the notation. There are now fully faithful 2-functors 
$$
\xymatrix{
\MDL(\K) = \Mnd(\Cmd(\K)) \ar[r] & \Mnd(\K^\two) \ar[r] & 
(\K^\two)^\two = \K^{\two\x\two} }
$$
sending an object $(K,t,c,\lambda)$ of $\MDL(\K)$ to the commutative square of
forgetful maps: 
$$
\xymatrix{
K^{(t,c)} \ar[r] \ar[d] & K^c \ar[d] \\
K^t \ar[r] & K }
$$
If a mixed distributive law in the 2-category $\Cat$ of Categories; Functors;
Natural Transformations is induced by an entwining structure in $\Bim$, then
the associated category $K^{(t,c)}$ is known as the category of `entwined
modules' \cite{BrzMaj}.   

The situation of distributive laws between two monads is not completely
analogous to the mixed case above. If $\cK$ admits Eilen\-berg-Moore
constructions for monads, then there is still a fully faithful 2-functor from
$\Mnd(\Mnd(\cK))$ -- considered as the 2-category of distributive laws between
two monads -- to $\cK^{\two\x\two}$.  The main difference in this
``non-mixed'' case is that while one of the monads lifts to the
Eilenberg-Moore object of the other, the other extends to the Kleisli object
of the first.  

In order to treat algebra extensions by weak bialgebras \cite{BNSz}, entwining
structures were generalized in \cite{CaeDeG} to `weak entwining structures',
which are better called `weak mixed distributive laws' if working in general
2-categories. A weak mixed distributive law  in a 2-category $\cK$ also
consists of a monad $(K,t)$ and a comonad $(K,c)$, together with a 2-cell $tc
\1c {} ct$, but the compatibility axioms with the unit of the monad and the
counit of the comonad are weakened. The corresponding notion of weak
distributive law between two monads is discussed in \cite{Str:weak_dl}. The
aim of this paper is to extend to weak (mixed) distributive laws the
standard results for ordinary (mixed) distributive laws sketched above. 

We are not aware of any characterization of a (mixed) weak distributive law as
a monad (or as a comonad) in some 2-category. Instead, in this note we observe
that a mixed weak distributive law in an arbitrary 2-category $\cK$ can be
described as a compatible pair consisting of a comonad in a 2-category
$\Mnd^\iota(\cK)$, extending $\Mnd(\cK)$, and a monad in
$\Cmd^\pi(\cK):=\Mnd^\iota(\cK_*)_*$, cf. \cite{weak_EM}. This observation is 
used in Section \ref{sec:wentw}  to define a 2-category $\WDL(\cK)$, whose
0-cells are weak mixed distributive laws in $\cK$ and whose 1-cells and
2-cells are also compatible pairs of 1-cells and 2-cells, respectively, in
$\Mnd(\Cmd^\pi(\cK))$ and $\Cmd(\Mnd^\iota(\cK))$. By construction, the
2-category $\WDL(\cK)$ comes equipped with 2-functors $\WDL(\cK) \to
\Cmd(\Mnd^\iota(\cK))$ and $\WDL(\cK) \to \Mnd(\Cmd^\pi(\cK))$, and indeed
$\WDL(\K)$ can be seen as a sort of ``intersection'' of $\Mnd(\Cmd^\pi(\K))$
and $\Cmd(\Mnd^\iota(\K))$.  

Although the 2-categories $\Mnd^\iota(\K)$ and $\Cmd^\pi(\K)$ do not embed in
$\K^\two$, they embed in a sort of ``weak'' version of $\K^\two$,
corresponding to the ``weak liftings'' studied in \cite{weak_EM} and 
\cite{Steve_talk}. Perhaps surprisingly, the 2-category $\WDL(\K)$ does
still embed in $\K^{\two\x\two}$. We prove this, and characterize the image of
the embedding, as well as describing how this relates to the weak liftings
described above.  

If a weak mixed distributive law in $\Cat$ is induced by a weak entwining
structure between a $k$-algebra $t$ and a $k$-coalgebra $c$, then the four
objects occurring in its image in $\Cat^{\two\x \two}$ are the category of
$k$-modules, the category of $t$-modules, the category of $c$-comodules and
the category of so-called weak entwined modules  \cite{CaeDeG},
\cite{Brz:coring}. 
Important examples of weak entwining structures are associated with
Doi-Koppinen data over weak bialgebras \cite{Bohm:Doi}. The corresponding weak
entwined modules include various Hopf type modules over weak bialgebras -- such
as (relative) Hopf modules and Yetter-Drinfel'd modules -- so in particular
graded modules over groupoid graded algebras (cf.\cite{CaeDeG:wGal}). More
exotic weak distributive laws, behind which there are no Doi-Koppinen data,
were constructed in \cite{BrTuWr}.

By a formal dualization of the above results on weak mixed distributive laws,
one can also define a 2-category whose 0-cells are weak distributive laws
between two monads. If Eilen\-berg-Moore constructions for monads exist and
also idempotent 2-cells split in $\cK$, then we obtain a fully faithful
2-functor from it to $\cK^{\two\x \two}$.   
\medskip
 
{\bf Notation.}
We assume that the reader is familiar with the basic theory of
2-categories. For a review of the required notions (such as 2-categories,
2-functors, 2-adjunctions, monads, adjunctions and Eilen\-berg-Moore
construction in a 2-category) we refer to the article \cite{KeSt}. In a
2-category $\cK$, horizontal composition is denoted by juxtaposition and
vertical composition is denoted by a dot. 
We say that in $\cK$ {\em idempotent 2-cells split} provided that  for any
2-cell $\Theta:V \1c{} V$ such that $\Theta. \Theta = \Theta$, there exist a
1-cell ${\widehat V}$ and 2-cells $\pi:V \1c{} {\widehat V}$ and
$\iota:{\widehat V} \1c{} V$, such that $\pi . \iota ={\widehat V}$ and $\iota
. \pi =\Theta$.  
\medskip

{\bf Acknowledgement.}
GB would like to thank the organizers of the Conference in Hopf
algebras and noncommutative algebra in the honor of Mia Cohen, in Sde-Boker,
May 2010, for a generous invitation and an unforgettable first time in Israel. 
She also acknowledges financial support of the Hungarian Scientific Research
Fund OTKA, grant no. F67910.  
All authors are grateful for partial support from the Australian
Research Council, project DP0771252.


\section{The 2-category of mixed weak distributive laws}
\label{sec:wentw}

A mixed weak distributive law in a 2-category \K consists of a monad
$(t,\mu,\eta)$ and a comonad $(c,\delta,\varepsilon)$ on an object $K$ of \K,
along with a 2-cell $\lambda:tc\to ct$ making the following diagrams commute: 
\begin{eqnarray}
\xymatrix @C=25pt @R=19pt{
ttc \ar[r]^-{\mu c}
\ar[d]_-{t\lambda}&
tc \ar[dd]^-{\lambda}\\
tct \ar[d]_-{\lambda t}&\\
ctt\ar[r]^-{c\mu}&
ct}\qquad
\xymatrix  @C=25pt @R=19pt{
tc \ar[r]^-{t \delta} \ar[dd]_-\lambda&
tcc\ar[d]^-{\lambda c}\\
&ctc\ar[d]^-{c \lambda}\\
ct \ar[r]^-{\delta t}&
cct}\qquad
\xymatrix @R=9pt  @C=25pt {
c\ar[r]^-\delta \ar[ddd]_-{\eta c}&
cc\ar[d]^-{c \eta c}\\
&ctc \ar[d]^-{c\lambda}\\
&cct 
\ar[d]^-{c \varepsilon t}\\
tc \ar[r]^-\lambda&
ct}\qquad
\xymatrix @R=9pt  @C=25pt{
tc \ar[r]^-\lambda 
\ar[d]_-{t \eta c}&
ct \ar[ddd]^-{\varepsilon t}\\
ttc \ar[d]_-{t\lambda}&\\
tct \ar[d]_-{t \varepsilon t}&\\
tt \ar[r]^-\mu&t.}
\label{eq:wentw_m}\label{eq:wentw_d} \label{eq:wentw_u}\label{eq:wentw_e}
\end{eqnarray}

A 2-functor, or more generally a pseudofunctor, $F:\K\to\L$ sends weak
distributive laws in \K to weak distributive laws in \L. In particular, a
representable 2-functor $\K(X,-)$ sends a weak distributive law as above to a 
weak distributive law in \Cat on the category $\K(X,K)$.

The main difference between mixed distributive laws and weak ones is that 
$(c,\lambda)$ is no longer a morphism of monads from $(K,t)$ to $(K,t)$,
and so no longer induces a lifting $\overline{c}:K^t\to K^t$ when the
Eilenberg-Moore object $K^t$ exists. Nonetheless, we shall see that 
$(c,\lambda)$ does induce a {\em weak lifting}
$\overline{c}:K^t\to K^t$; equivalently, $(c,\lambda)$ is a {\em weak} morphism
of monads $(K,t)\to(K,t)$ \cite{weak_EM}. These weak morphisms of monads 
are the 1-cells of a 2-category $\Mnd^\iota(\K)$ whose objects are monads in \K.

\begin{theorem}[\cite{weak_EM}, Corollary 1.4]
\label{thm:Mnd^i} 
For any 2-category $\cK$, the following data constitute a 2-category, to be
denoted by $\Mnd^\iota(\cK)$.
\begin{itemize}
\item[] \underline{0-cells} are monads $(K,t)$ in $\cK$.
\item[] \underline{1-cells} $(K,t)\1c{} (K',t')$ are pairs, consisting of a
  1-cell $x:K\1c {} K'$ and a 2-cell $\xi:t' x \1c {}  x
  t$ in $\cK$ such that the first diagram in \eqref{eq:1-cell} commutes. 
\item[] \underline{2-cells} $(x,\xi)\1c{} (x',\xi')$ are 2-cells
$\omega:x\1c{}x'$ in $\cK$, rendering commutative the second diagram in
  \eqref{eq:1-cell}. 
\item[] \underline{Horizontal and vertical compositions} are the same as in
  $\cK$. 
\end{itemize}
The 2-category $\Mnd^\iota(\cK)$ contains $\Mnd(\cK)$ as a vertically full
2-subcategory. 
\begin{equation}\label{eq:1-cell}
\xymatrix{
t't'x\ar[r]^-{t'\xi}\ar[d]_-{\mu'x}&
t'xt\ar[r]^-{\xi t}&
xtt\ar[d]^-{x\mu}\\
t'x\ar[rr]^-{\xi}&&xt}\qquad
\xymatrix @C=5pt{
t'x\ar[rrr]^-\xi \ar[d]_-{t'\eta'x}&&&
xt\ar[rrr]^-{\omega t}&&&
x't\\
t't'x\ar[rr]^-{t'\xi}&&
t'xt\ar[rr]^-{t'\omega t}&&
t'x't\ar[rr]^-{\xi' t}&&
x'tt\ar[u]_-{x'\mu}}
\end{equation}
\end{theorem}
In \cite[Corollary 1.4]{weak_EM} another 2-category $\Mnd^\pi(\cK)$ was
introduced, with the same 0- and 1-cells as in $\Mnd^\iota(\cK)$ but different
2-cells. $\Mnd^\pi(\cK)$ also  contains $\Mnd(\cK)$ as a vertically full
2-subcategory. 

Similarly, for a weak mixed distributive law $(K,t,c,\lambda)$, we have only
a weak morphism of comonads $(t,\lambda):(K,c)\to(K,c)$, and only a weak
lifting of $t$ to $K^c$. There is a 2-category $\Cmd^\pi(\K)=\Mnd^\iota(\K_*)_*$
of comonads in \K and weak morphisms of comonads; once again, it contains
$\Cmd(\K)$ as a vertically full sub-2-category.

Our aim is to construct a 2-category of weak mixed distributive laws in any
2-category $\cK$; that is, a 2-category $\WDL(\cK)$ whose objects are weak
mixed distributive laws. Our starting point is the following lemma. 

\begin{lemma}\label{lem:entw_0-cells}
For a monad $(K,t)$, a comonad $(K,c)$ and a 2-cell $\lambda:tc\1c{} ct$ in any
2-category $\cK$, the following statements are equivalent.
\begin{itemize}
\item[{(i)}] $(K,t,c,\lambda)$ is a weak mixed distributive law;
\item[{(ii)}] $((K,t),(c,\lambda))$ is a comonad in $\Mnd^\iota(\cK)$ and
  $((K,c),(t,\lambda))$ is a monad in  $\Cmd^\pi(\cK)$. 
\end{itemize}
\end{lemma}

\begin{proof}
The first axiom in \eqref{eq:wentw_m} expresses the requirement
that $(c,\lambda):(K,t)\1c{}(K,t)$ is a 1-cell in $\Mnd^\iota(\cK)$ and the
second axiom in \eqref{eq:wentw_d} means that $(t,\lambda):(K,c)
\1c{} (K,c)$ is a 1-cell in $\Cmd^\pi(\cK)$. The third axiom in
\eqref{eq:wentw_u} means that $\eta:1_{(K,c)}\1c {} (t,\lambda)$ is a 2-cell
in $\Cmd^\pi(\cK)$ and the last axiom in \eqref{eq:wentw_e} means
that $\varepsilon:(c,\lambda)\1c {} 1_{(K,t)}$ is a 2-cell in
$\Mnd^\iota(\cK)$. If these four conditions hold, then also
$\mu:(t,\lambda)(t,\lambda)\1c {} (t,\lambda)$ is a 2-cell in
$\Cmd^\pi(\cK)$. That is, the following diagram commutes.  
$$
\xymatrix @R=8pt{
t^2c\ar[r]^-{t^2\delta}\ar[d]_-{\mu c}&
t^2c^2\ar[r]^-{t\lambda c}\ar[d]^-{\mu c^2}&
(tc)^2\ar[r]^-{\lambda tc}&
ct^2c\ar[d]^-{c\mu c}\\
tc\ar[r]^-{t\delta}\ar[d]_-\lambda&
tc^2\ar[rr]^-{\lambda c}&&
ctc\ar[d]^-{c\lambda}\\
ct\ar[rrr]^-{\delta t}\ar@{=}[rrrd]&&&
c^2t\ar[d]^-{c\varepsilon t}\\
&&&ct}
$$
The top right region commutes by the first axiom, and the region
below it commutes by the second axiom in \eqref{eq:wentw_e}. The
top left square commutes by naturality and the triangle commutes
by a counitality axiom of a comonad. Symmetrically, \eqref{eq:wentw_m} implies
that $\delta:(c,\lambda)\1c {} (c,\lambda)(c,\lambda)$ is a 2-cell in
$\Mnd^\iota(\cK)$.  
\end{proof}

The situation in Lemma \ref{lem:entw_0-cells} (ii) is distinguished among the
other possibilities in the following sense. 
\begin{lemma}
For a monad $(K,t)$, a comonad $(K,c)$, and a 2-cell $\lambda:tc \1c {} ct$ in
any 2-category $\cK$, the following assertions are equivalent.
\begin{itemize}
\item[{(i)}] $(K,t,c,\lambda)$ is mixed distributive law; 
\item[{(ii)}] $((K,t),(c,\lambda))$ is a comonad in $\Mnd(\cK)$;
\item[{(iii)}] $((K,c),(t,\lambda))$ is a monad in $\Cmd(\cK)$;
\item[{(iv)}] $((K,c),(t,\lambda))$ is a monad in $\Cmd^\iota(\cK)$ and
  $((K,t),(c,\lambda))$ is a comonad in  $\Mnd^\pi(\cK)$;
\item[{(v)}] $((K,c),(t,\lambda))$ is a monad in $\Cmd^\pi(\cK)$ and
  $((K,t),(c,\lambda))$ is a comonad in $\Mnd^\pi(\cK)$;
\item[{(vi)}] $((K,c),(t,\lambda))$ is a monad in $\Cmd^\iota(\cK)$ and
  $((K,t),(c,\lambda))$ is a comonad in  $\Mnd^\iota(\cK)$.
\end{itemize}
\end{lemma}
\begin{proof}
Equivalence of (i), (ii) and (iii) is well-known, see e.g. \cite[Proposition
6.3 and Corollary 6.6]{PowWat}.  
Assertions (ii) and (iii) trivially imply any of (iv), (v) and (vi).
The counit $\varepsilon$ of $(K,c)$ is a 2-cell $(c,\lambda)\1c{}
1_{(K,t)}$ in $\Mnd^\pi(\cK)$ if and only if $\varepsilon t . \lambda =
t\varepsilon$ and the unit $\eta$ of $(K,t)$ is a 2-cell $1_{(K,c)}
\1c{} (t,\lambda)$ in $\Cmd^\iota(\cK)$ if and only if it 
$\lambda .\eta c= c\eta$. Hence (iv)$\Rightarrow$(i).   
If (v) holds then $\lambda$ obeys the third axiom in \eqref{eq:wentw_u} and
$\varepsilon t . \lambda = t\varepsilon$, hence by counitality of $\delta$
also $\lambda.\eta c= c\eta$. This proves (v)$\Rightarrow$(i) and
(vi)$\Rightarrow$(i) follows symmetrically.   
\end{proof}

The next two lemmas are preparatory to our definition of 1-cells and 2-cells
in the 2-category $\WDL(\K)$. 

\begin{lemma}\label{lem:entw_1-cells}
For weak mixed distributive laws $(K,t,c,\lambda)$ and $(K',t',c',\lambda')$
in a 2-category $\cK$, consider a 1-cell $(x,\xi):(K,t)\1c{} (K',t')$ in
$\Mnd(\cK)$ and a 1-cell $(x,\zeta):(K,c)\1c{} (K',c')$ in
$\Cmd(\cK)$. The following are equivalent. 
\begin{itemize}
\item[{(i)}] $\xi:(t',\lambda')(x,\zeta)\1c{} (x,\zeta)(t,\lambda)$ is a 2-cell
  in $\Cmd^\pi(\cK)$; that is, the following diagram commutes;
  \begin{equation}\label{eq:wdl_morphism}
\xymatrix{
t'xc^2 \ar[r]^{t'\zeta c} & 
t'c'xc \ar[r]^{\lambda' xc} & 
c't'xc \ar[r]^{c'\xi c} & 
c'xtc \ar[r]^{c'x\lambda} & 
c'xct \ar[d]^{c'x\varepsilon t} \\
t'xc \ar[u]^{t'x\delta} \ar[r]^-{\xi c} & 
xtc \ar[rr]^-{x\lambda} && 
xct \ar[r]^-{\zeta t} & c'xt}
  \end{equation}
\item[{(ii)}] $\zeta:(x,\xi)(c,\lambda)\1c{} (c',\lambda')(x,\xi)$ is a 2-cell 
  in  $\Mnd^\iota(\cK)$;  that is, the following diagram commutes.
  \begin{equation}
\xymatrix{
t'xtc \ar[r]^{t'x\lambda} & 
t'xct \ar[r]^{t'\zeta t} & 
t'c'xt \ar[r]^{\lambda' xt} & 
c't'xt \ar[r]^{c'\xi t} & 
c'xtt \ar[d]^{c'x\mu} \\
t'xc \ar[u]^{t'x\eta c} \ar[r]^-{\xi c} & 
xtc \ar[rr]^-{x\lambda} && 
xct \ar[r]^-{\zeta t} & c'xt  
}
 \end{equation}
When these conditions hold, we say that $(x,\xi,\zeta)$ is a morphism of weak
distributive laws from $(K,t,c,\lambda)$ to $(K',t',c',\lambda')$.  
\end{itemize}
\end{lemma}

\begin{proof}
This follows by commutativity of 
$$\xymatrix @R=15pt{
t'xc  \ar[r]^-{t'x\delta} \ar[ddd]_-{t'x\eta c} &
t'xc^2 \ar[d]^-{t'xc\eta c} \ar[r]^-{t'\zeta c} & 
t'c'xc \ar[r]^-{\lambda'xc}&
c't'xc \ar[r]^-{c'\xi c}&
c'xtc \ar[d]^-{c'xt\eta c} \ar[r]^-{c'x\lambda}&
c'xct \ar[ddd]^-{c'x\varepsilon t}\\
&t'xctc \ar[d]^-{t'xc\lambda} &&&
c'xt^2 c \ar[d]^{c'xt\lambda} & \\
&t'xc^2 t \ar[d]^-{t'xc\varepsilon t}&&&
c'xtct \ar[d]^-{c'xt\varepsilon t}&\\
t'xtc \ar[r]^{t'x\lambda} &
t'xct \ar[r]^{t'\zeta t} &
t'c'xt \ar[r]^{\lambda'xt} &
c't'xt \ar[r]^{c'\xi t} &
c'xt^2 \ar[r]^{c'x\mu} &
c'xt }$$
in which the large central region commutes by naturality, and the other two
regions by the weak distributive law axioms.
\end{proof}

\begin{lemma}
For morphisms
$(x,\xi,\zeta),(x',\xi',\zeta'):(K,t,c,\lambda)\to(K',t',c',\lambda')$ of weak 
distributive laws  in a 2-category $\cK$, and a 2-cell $\omega:x\to x'$, the
following conditions are equivalent.
\begin{itemize}
\item[{(i)}] $\omega:(x,\xi)\to(x',\xi')$ is a 2-cell in $\Mnd(\K)$ and 
$\omega:(x,\zeta)\to(x,\zeta')$ is a 2-cell in $\Cmd(\K)$;
\item[{(ii)}] $\omega:((x,\xi),\zeta)\to((x',\xi'),\zeta')$ is a 2-cell in 
$\Cmd(\Mnd^\iota(\K))$ and $\omega:((x,\zeta),\xi)\to((x',\zeta'),\xi')$
is a 2-cell in $\Mnd(\Cmd^\pi(\K))$.
\end{itemize}
We then say that $\omega$ is a 2-cell $(x,\xi,\zeta)\to(x',\xi',\zeta')$ of
weak distributive laws. 
\end{lemma}

\begin{proof}
In each case the conditions assert the commutativity of the squares
$$\xymatrix{
t'x \ar[r]^{t'\omega} \ar[d]_{\xi} & 
t'x' \ar[d]^{\xi'} \\
xt \ar[r]^-{\omega t} & x't }\qquad
\xymatrix{
xc \ar[r]^{\omega c} \ar[d]_{\zeta} & 
x'c \ar[d]^{\zeta'} \\
c'x \ar[r]^-{c'\omega} & c'x'. }
$$
\end{proof}
We now deduce: 
\begin{theorem}\label{thm:wentw}
For any 2-category $\cK$, the weak distributive laws in \K along with the 
morphisms and 2-cells defined above give a 2-category $\WDL(\K)$ with
composition performed as in $\Mnd(\K)$ and $\Cmd(\K)$.
\end{theorem}

\section{A fully faithful embedding for weak mixed distributive laws}

For any 2-category \K, there is a fully faithful 2-functor $Y:\K\to\WDL(\K)$
which equips any object $X$ with the identity monad, the identity comonad, and
the identity distributive law between them; we write $YX=(X,1,1,1)$. 

If $(K,t,c,\lambda)$ is any weak distributive law, a morphism in $\WDL(\K)$ 
from $YX$ to $(K,t,c,\lambda)$ consists of a morphism $a:X\to K$ in \K, equipped
with an action $\alpha:ta\to a$ of the monad $t$ and a  coaction $\gamma:a\to
ca$ of the comonad $c$, satisfying the compatibility condition asserting that
the diagram  
$$
\xymatrix{
tca\ar[rr]^-{\lambda a}&&cta\ar[d]^-{c\alpha}\\
ta\ar[u]^-{t\gamma}\ar[r]^-\alpha&a\ar[r]^-\gamma&ca
}
$$
commutes. We then say that $(a,\alpha,\gamma)$ is a mixed
$(K,t,c,\lambda)$-algebra with domain $X$. A morphism of mixed
$(K,t,c,\lambda)$-algebras is a 2-cell in \K compatible with the action and
coaction. Thus the category of mixed algebras with domain $X$ is the
hom-category  
$$
\WDL(\K)(YX,(K, t,c,\lambda)).
$$

In particular, we may take $\K=\Cat$ and $X=1$, the terminal category. Then
$a$ is just an object of $\cK$, and the action and coaction amount to a
$t$-algebra structure and a $c$-coalgebra structure in the usual sense, while
the compatibility condition has exactly the form of the diagram displayed
above. We then write $K^{(t,c)}$ for the category of mixed
$(K,t,c,\lambda)$-algebras.   

In fact we can also recover the general notion of mixed algebra from this
particular one: for a weak distributive law $(K, t,c,\lambda)$ in a 2-category
\K, applying the representable 2-functor $\K(X,-):\K\to\Cat$ gives a weak
distributive law   
$$
(\K(X,K),\K(X,t),\K(X,c),\K(X,\lambda))
$$ 
in \Cat whose category of mixed algebras is just the category of mixed
$(K,t,c,\lambda)$-algebras with domain $X$. This passage from $X$ to the
category of mixed $(K,t,c,\lambda)$-algebras with domain $X$ defines a
2-functor $\K\op\to\Cat$, and if this 2-functor is representable, say as  
$$
\K(X,K)^{(\K(X,t),\K(X,c))}\cong \K(X,K^{(t,c)}),
$$
we call the representing object $K^{(t,c)}$ the mixed Eilenberg-Moore object.
(Clearly, if $\cK=\Cat$, we re-cover the above category $K^{(t,c)}$ of
mixed $(K,t,c,\lambda)$-algebras with domain $1$.)  

\begin{theorem}
\K has mixed Eilenberg-Moore objects if and only if $Y:\K\to\WDL(\K)$ has a
right 2-adjoint. 
\end{theorem}

\proof
This follows from the isomorphism 
$$
\WDL(\K)(YX,(K, t,c,\lambda))\cong \K(X,K)^{(\K(X,t),\K(X,c))}
$$
since the mixed Eilenberg-Moore object is defined as a representing object for
the right hand side, while the right adjoint is defined as a representing
object for the left hand side. 
\endproof

\begin{theorem}
Suppose that \K has Eilenberg-Moore objects for monads and comonads, and that 
idempotent 2-cells split. Then \K has mixed Eilenberg-Moore objects.
\end{theorem}

\proof
Form the Eilenberg-Moore object $K^c$. Since $(t,\lambda)$ is a monad in 
$\Cmd^\pi(\K)$, it has a weak lifting to a monad $\overline{t}$ on $K^c$ 
cf. \cite[Proposition 5.7]{weak_EM}. The Eilenberg-Moore object
$(K^c)^{\overline{t}}$ is the desired mixed Eilenberg-Moore object.
Symmetrically, the comonad $(c,\lambda)$ has a weak lifting to a
monad $\overline{c}$ on $K^t$, and $(K^t)^{\overline{c}}$ also gives the mixed 
Eilenberg-Moore object. 
\endproof

\begin{theorem}
Suppose that \K has Eilenberg-Moore objects for monads and comonads and
that idempotent 2-cells split. Then $\WDL(\K)$ has a fully faithful embedding
into $\K^{\two\x\two}$. 
\end{theorem}

\proof
The embedding will send an object $(K, t,c,\lambda)$ to the square
$$
\xymatrix{
K^{(t,c)} \ar[r]^{\overline{v}} \ar[d]_{\overline{u}} & K^t \ar[d]^{u} \\
K^c \ar[r]^-{v} & K }
$$
of Eilenberg-Moore objects and forgetful morphisms. 
In order to conclude that it extends to the stated embedding, we need to
establish isomorphisms between the respective hom categories of $\WDL(\K)$ and 
$\K^{\two\x\two}$. 

A 1-cell in $\WDL(\K)$ consists of a compatible pair of a monad morphism
and a comonad morphism.   
If $(K',t',c',\lambda')$ is another object of $\WDL(\K)$, then consider a
morphism $x:K\to K'$. To give a monad morphism $(x,\xi):(K,t)\to(K',t')$ is
equivalently to give a lifting $x^\xi:K^t\to K'^{t'}$ of $x$; this in turn is
equivalent to the fact that for any $t$-algebra with domain $X$, consisting of
a morphism $a:X\to K$ with an action $\alpha:ta\to a$, the composite  
$$
\xymatrix{
t'xa \ar[r]^{\xi a} & xta \ar[r]^{x\alpha} & xa }
$$
makes $xa$ into a $t'$-algebra with domain $X$. Similarly, to give a comonad
morphism $(x,\zeta):(K,c)\to(K',c')$ is equivalently to give a lifting
$x^\zeta:K^c\to K'^{c'}$ of $x$; this is equivalent to the fact that for any
$c$-coalgebra $(a,\gamma)$ with domain $X$,  the composite
$$
\xymatrix{
xa \ar[r]^{x\gamma} & xca \ar[r]^{\zeta a} & c'xa }
$$
makes $xa$ into a $c'$-coalgebra with domain $X$.

We should check that $(x,\xi,\zeta)$ is a morphism in $\WDL(\K)$ if and only
if $x^\xi$ and $x^\zeta$ have a common lifting $x^{(\xi,\zeta)}:K^{(t,c)}\to
K'^{(t',c')}$ (giving rise to a 1-cell in $\cK^{2\x2}$);   
in other words, if and only for any mixed $(K,t,c,\lambda)$-algebra
$(a,\alpha,\gamma)$ with domain $X$, the induced $t'$-algebra and
$c'$-coalgebra structures on $xa$ together give a
$(K',t',c',\lambda')$-algebra. But we can regard $(a,\alpha,\gamma)$ as a 
morphism in  $\WDL(\K)$ from $YX$ to $(K,t,c,\lambda)$, and now composing with
$(x,\xi,\zeta)$ gives a morphism $YX\to(K',t',c',\lambda')$ which gives the
desired $(K',t',c',\lambda')$-algebra.  

Conversely, suppose that $x^\xi$ and $x^\zeta$ have a common lifting
$x^{(\xi,\zeta)}$. The composite  
$$
\xymatrix{
tc \ar[r]^{t\delta} & tc^2 \ar[r]^{\lambda c} & ctc \ar[r]^{\varepsilon tc} & tc}
$$
is idempotent, and splits to give a morphism $a:K\to K$ equipped with an
epimorphism $\pi:tc\to a$ and monomorphism $\sigma:a\to tc$. Then $a$ has a
mixed $(K,t,c,\lambda)$-algebra structure $(a,\alpha,\gamma)$ coming from  
\begin{equation}\label{eq:alpha_gamma}
\xymatrix @R0pc {
ta \ar[r]^{t\sigma} &t^2 c  \ar[r]^{\mu c} & tc \ar[r]^{\pi} & a} 
\quad \textrm{and}\quad
\xymatrix @R0pc {
a \ar[r]^{\sigma} & tc \ar[r]^{t\delta} & tc^2 \ar[r]^{\lambda c} & ctc
\ar[r]^{c\pi} &  ca. }
\end{equation} 
Applying the lifting $x^{(\xi,\zeta)}$ gives a mixed
$(K',t',c',\lambda')$-algebra structure on $xa$ via the maps 
$$
\xymatrix @R0pc {
t'xa \ar[r]^{\xi a} & xta \ar[r]^{x\alpha} & xa}
\quad \textrm{and}\quad
\xymatrix @R0pc {
xa \ar[r]^{x\gamma} & xca \ar[r]^{\zeta a} & c'xa. }
$$
Write $\eta_1:c\to a$ and $\varepsilon_1:a\to t$ for the composites
$$
\xymatrix @R0pc {
c \ar[r]^{\eta c} & tc \ar[r]^{\pi} & a & {\textnormal{and}} &
a \ar[r]^{\sigma} & tc \ar[r]^{t\varepsilon} & t. }
$$
In view of \eqref{eq:wdl_morphism}, 
$(x,\xi,\zeta)$ will be a morphism in $\WDL(\K)$ provided that the exterior of  
$$
\xymatrix{
t'xc^2 \ar[rd]^{t'xc\eta_1} \ar[rr]^{t'\zeta c} && t'c'xc \ar[r]^{\lambda' xc} & 
c't'xc \ar[rr]^{c'\xi c} && c'xtc \ar[dl]^{c'xt\eta_1} \ar[d]^{c'x\lambda} \\
& t'xca \ar[r]^{t'\zeta a} & t'c'xa \ar[r]^{\lambda' xa} & c't'xa
\ar[r]^{c'\xi a} & c'xta \ar[d]_{c'x\alpha} & 
c'xct \ar[d]^{c'x\varepsilon t} \\
t'xc \ar[r]^{t'x\eta_1} \ar[uu]^{t'x\delta} \ar[dd]_{\xi c} & 
t'xa \ar[u]_{t'x\gamma} \ar[d]^{\xi a} &&& c'xa \ar[r]^{c'x\varepsilon_1} & c'xt \\
& xta \ar[rr]^{x\alpha} \ar[dr]^{xt\gamma}  && xa \ar[r]^{x\gamma} & 
xca \ar[u]^{\zeta a}  \ar[rd]^{xc\varepsilon_1} \\
xtc \ar[ur]_{xt\eta_1} \ar@/_2pc/[rrrrr]_{x\lambda}  && xtca \ar[r]^{x\lambda
  a} & xcta \ar[ur]^{xc\alpha} && xct \ar[uu]_{\zeta t}  
}
$$
commutes. The central region commutes since $xa$ is a mixed algebra,  the
region just below it because $a$ is a mixed algebra, the top region and the
lower corners by naturality. Thus we just need to show that the lower region
and the two upper corners commute. These will follow from commutativity of the
following three diagrams: 
\begin{equation}\label{eq:lemma}
\xymatrix{
ta \ar[r]^{t\gamma} & tca \ar[r]^{\lambda a} & cta \ar[r]^{c\alpha} & ca
\ar[d]^{c\varepsilon_1} \\ 
tc \ar[u]^{t\eta_1} \ar[rrr]^-{\lambda} &&& ct }
\xymatrix{ 
c \ar[r]^{\eta_1} \ar[d]_{\delta} & a \ar[d]^{\gamma} \\
c^2\ar[r]^-{c\eta_1} & ca}
\xymatrix{
tc \ar[r]^{t\eta_1} \ar[d]_{\lambda} & ta \ar[r]^{\alpha} & a \ar[d]^{\varepsilon_1} \\
ct \ar[rr]^-{\varepsilon t} && t }
\end{equation}
which involve a single weak distributive law ($x$ does not appear). We prove
them as a separate lemma below. 

We now turn to fully faithfulness on 2-cells. Let $(x,\xi,\zeta)$ and
$(x',\xi',\zeta')$ be 1-cells from $(K,t,c,\lambda)$ to $(K',t',c',\lambda')$,
with induced liftings $x^\xi$, $x'^{\xi'}$, and so on. A 2-cell $\omega:x\to
x'$ in \K lifts to a 2-cell $\omega^\xi:x^\xi\to x'^{\xi'}$ if and only if
$\omega$ is a monad 2-cell $(x,\xi)\to(x',\xi')$; similarly it lifts to a
2-cell $\omega^\zeta:x^\zeta\to x'^{\zeta'}$ if and only if it is a comonad
2-cell $(x,\zeta)\to(x',\zeta')$. Suppose both of these conditions hold. There
is no further condition for $\omega$ to be a 2-cell
$(x,\xi,\zeta)\to(x',\xi',\zeta')$; thus it remains to show that $\omega^\xi$
and $\omega^\zeta$ have a  common lifting
$\omega^{(\xi,\zeta)}:x^{(\xi,\zeta)}\to x'^{(\xi',\zeta')}$
(giving rise to a 2-cell in $\cK^{2\x 2}$). 
For this, it will suffice to show that for any mixed $(K,t,c,\lambda)$-algebra
$(a,\alpha,\gamma)$, if we form the induced mixed
$(K',t',c',\lambda')$-algebra structures on $xa$ and $x'a$, then $\omega
a:xa\to x'a$ will be a morphism of mixed $(K',t',c',\lambda')$-algebras. This
follows by regarding $(a,\alpha,\gamma)$ as a morphism $YX\to(K,t,c,\lambda)$
in $\WDL(\K)$, and then composing this with the 2-cell
$\omega:(x,\xi,\zeta)\to(x',\xi',\zeta')$. This gives a 2-cell
$(x,\xi,\zeta)(a,\alpha,\gamma)\to(x',\xi',\zeta')(a,\alpha,\gamma)$ with
domain $YX$ and codomain $(K',t',c',\lambda')$, and this is the desired
morphism of $(K',t',c',\lambda')$-algebras.  
\endproof

\begin{lemma} For any weak mixed distributive law $(K,t,c,\lambda)$ in
  a 2-category $\cK$ with split idempotent 2-cells, the diagrams in
  \eqref{eq:lemma} commute.  
\end{lemma}

\begin{proof}
In the diagrams
$$
\xymatrix @R=15pt{
tc \ar[r]^-{t\delta}\ar[dd]_-{t\delta}&
tc^2\ar[r]^-{\lambda c} \ar[d]^-{tc\delta}&
ctc\ar[r]^-{\varepsilon tc}\ar[d]^-{ct\delta}&
tc\ar[d]^-{t\delta}\\
&tc^3\ar[r]^-{\lambda c^2}\ar@{}[rd]|(.4){(*)}&
ctc^2\ar[r]^-{\varepsilon tc^2}\ar[d]^-{c\lambda c}&
tc^2\ar[d]^-{\lambda c}\\
tc^2\ar[ru]^-{t\delta c}\ar[r]^-{\lambda c}&
ctc\ar[r]^-{\delta tc}\ar@{=}@/_1.2pc/[rr]&
c^2tc\ar[r]^-{\varepsilon ctc}&
ctc
}\qquad
\xymatrix  @R=15pt{
tc \ar[r]^-{t\delta}\ar[d]_-{t\delta}&
tc^2\ar[r]^-{\lambda c}\ar[d]^-{t\delta c}\ar@{}[rdd]|{(*)}&
ctc\ar[dd]^-{\delta tc}\ar@{=}[rdd]&\\
tc^2\ar[r]^-{tc\delta}\ar[d]_-{\lambda c}&
tc^3\ar[d]^-{\lambda c^2}&&\\
ctc\ar[r]^-{ct\delta }&
ctc^2\ar[r]^-{c\lambda c}&
c^2tc \ar[r]^-{c\varepsilon tc}&
ctc}
$$
\ 

\noindent
the marked regions commute by the second  axiom  of a weak mixed distributive
law. The other regions commute by naturality, coassociativity and counitality
of the comonad $c$. This proves that $\lambda c. t\delta.\sigma.\pi=\lambda
c. t\delta = c\sigma.c\pi.\lambda c. t\delta$. Hence recalling from
\eqref{eq:alpha_gamma} the formula of $\gamma$, and composing the top-right
path of the second diagram in \eqref{eq:lemma} by the monomorphism $c\sigma$
at the end, we obtain the morphism 
\begin{equation}\label{eq:diag2}
\big(\xymatrix{
c\ar[r]^-{\eta c}&tc \ar[r]^-{t\delta}&tc^2\ar[r]^-{\lambda c}&ctc}
\big)=\big(
\xymatrix{
c\ar[r]^-\delta&c^2\ar[r]^-{\eta c^2}&tc^2\ar[r]^-{\lambda c}&ctc}
\big).
\end{equation}
In the following diagram the square on the right commutes by the third axiom
of a weak mixed distributive law; the top left square commutes by
coassociativity of $\delta$ and the bottom left square commutes by naturality. 
$$
\xymatrix  @R=15pt{
c\ar[r]^-{\delta}\ar[d]_-{\delta}&c^2\ar[rr]^-{\eta c^2}\ar[d]^-{\delta c}&&
tc^2\ar[dd]^-{\lambda c}\\
c^2\ar[r]^-{c\delta}\ar[d]_-{c\eta c}&
c^3\ar[d]^-{c\eta c^2}&\\
ctc\ar[r]^-{ct\delta}&
ctc^2\ar[r]^-{c\lambda c}&
c^2tc\ar[r]^-{c\varepsilon tc}&
ctc
}
$$
Hence composing by $c\sigma$ the left-bottom path of the second diagram in
\eqref{eq:lemma}, we get the same morphism \eqref{eq:diag2}. This proves
commutativity of the second diagram in \eqref{eq:lemma}. 

The diagrams
$$
\xymatrix @C=15pt{
tc\ar[r]^-{t\delta}\ar[d]_-{t\eta c}&
tc^2\ar[rrr]^-{\lambda c}\ar[d]^-{t\eta c^2}&&&
ctc\ar[d]^-{\varepsilon tc}\\
t^2c\ar[r]^-{t^2\delta}&
t^2c^2\ar[r]^-{t\lambda c}&
(tc)^2\ar[r]^-{t\varepsilon tc}&
t^2c\ar[r]^-{\mu c}&
tc
}\qquad 
\xymatrix{
tc\ar[r]^-{t\delta}\ar@{=}[rd]&
tc^2\ar[r]^-{\lambda c}\ar[d]^-{tc\varepsilon}&
ctc\ar[r]^-{\varepsilon tc}\ar[d]^-{ct\varepsilon}&
tc\ar[d]^-{t\varepsilon}\\
&tc \ar[r]^-{\lambda}&
ct \ar[r]^-{\varepsilon t}&
t
}
$$
commute  by the last  axiom  of a weak mixed distributive law; naturality and
counitality of the comonad $c$. The equal paths around the first diagram
describe an idempotent morphism. Thus recalling the expression of $\alpha$
from \eqref{eq:alpha_gamma}, the third diagram in \eqref{eq:lemma} is seen to
be commutative.  

Finally, the first diagram in \eqref{eq:lemma} is identical to the outer
square in  
$$
\xymatrix @R=13pt {
ta\ar[r]^-{t\gamma}&
tca\ar[r]^-{\lambda a}&
cta\ar[r]^-{c\alpha}&
ca\ar[dd]^-{c\varepsilon_1}\\
tc\ar[u]^-{t\eta_1}\ar[r]^-{t\delta}\ar[d]_-{\lambda}&
tc^2\ar[u]_-{tc\eta_1}\ar[r]^-{\lambda c}&
ctc\ar[u]_-{ct\eta_1}\ar[d]^-{c\lambda}&&\\
ct\ar[rr]^-{\delta t}\ar@{=}@/_1.2pc/[rrr]&&
c^2t\ar[r]^-{c\varepsilon t}&
ct
}
$$
\ 

\noindent
Its top left and right regions commute by the second and third diagrams in
\eqref{eq:lemma}, respectively. The pentagon at the bottom left commutes by
the second  axiom  of a weak mixed distributive law. The middle square at the
top commutes by naturality and the bottom path is the identity morphism by the
counitality of $c$. 
\end{proof}

\section{A characterization of the image}

For the rest of the paper, we suppose that \K is a 2-category with
Eilenberg-Moore objects for monads and comonads, and splittings for idempotent
2-cells. We saw in the previous section that there is a fully
faithful 2-functor $\WDL(\K)\to\K^{\two\x\two}$; in this section we
characterize its image. 

First recall that for an adjunction $f\dashv u$, as well as the bijection
between 2-cells $fx\to y$ and $x\to uy$, there is also a bijection between
2-cells $mu\to n$ and $m\to nf$. Combining these two facts, we see that if
also $\overline{f}\dashv\overline{u}$ is an adjunction, then there are
bijections between 2-cells $fx\to y\overline{f}$, 2-cells $x\to
uy\overline{f}$, and 2-cells $x\overline{u}\to uy$. A 2-cell $fx\to
y\overline{f}$ and the corresponding 2-cell $x\overline{u}\to uy$ are said to
be {\em mates} \cite{KeSt}. Explicitly, the mates of $\phi:x\overline{u}\to
uy$ and $\psi:fx\to y\overline{f}$  are the composites
$$
\xymatrix @R0pc {
fx \ar[r]^{fx\overline{\eta}} & 
fx\overline{u}\overline{f} \ar[r]^{f\phi\overline{f}} & 
fuy\overline{f} \ar[r]^{\varepsilon y\overline{f}} & y\overline{f} }\qquad
\xymatrix @R0pc {
x\overline{u} \ar[r]^{\eta x\overline{u}} & 
ufx\overline{u} \ar[r]^{u \psi\overline{u}} & 
uy\overline{f}\overline{u} \ar[r]^{uy\overline{\varepsilon}} & uy}
$$
where as usual $\eta$ and $\varepsilon$ (possibly with a bar) denote the unit and
counit of an adjunction.

Consider a commutative square
$$
\xymatrix{
P \ar[r]^{\overline{v}} \ar[d]_{\overline{u}} & L \ar[d]^{u} \\
M \ar[r]^-{v} & K }
$$
in \K, and suppose that we have adjunctions $f\dashv u$,
$\overline{f}\dashv\overline{u}$, $v\dashv g$, and $\overline{v}\dashv
\overline{g}$.   

The identity 2-cell $v\overline{u}=u\overline{v}$ has mates
$\mone:fv\to\overline{v}\overline{f}$ and $\mtwo:\overline{u}\,\overline{g}\to
gu$ given by the following composites. 
$$
\xymatrix @R0pc {
fv \ar[r]^-{fv\overline{\eta}} & fv\overline{u}\overline{f} \ar@{=}[r] & 
fu\overline{v}\overline{f} \ar[r]^-{\varepsilon\overline{v}\overline{f}} & 
\overline{v}\overline{f}}\qquad 
\xymatrix @R0pc {
\overline{u}\,\overline{g} \ar[r]^-{\eta\overline{u}\,\overline{g}} &
gv\overline{u}\,\overline{g} \ar@{=}[r] & gu\overline{v}\,\overline{g}
\ar[r]^-{gu\overline{\varepsilon}} & 
gu }
$$
We shall use the following easy lemma:

\begin{lemma}
Let $\ell\dashv r$ and $\ell'\dashv r'$ be adjunctions, and
$\alpha:\ell\to\ell'$ a 2-cell with mate $\beta:r'\to r$. Then
$\check{\alpha}:\ell'\to\ell$ satisfies $\alpha\check{\alpha}=1$ if and only
if the mate $\hat{\beta}$ of $\check{\alpha}$ satisfies
$\hat{\beta}\beta=1$. \endproof 
\end{lemma}

\begin{proposition}
There is a bijection between 2-cells $\monesec:\overline{v}\overline{f}\to fv$
for which $\mone\monesec=1$ and 2-cells
$\mtworetr:gu\to\overline{u}\,\overline{g}$ for which  $\mtworetr\mtwo=1$.  
\end{proposition}

\proof
The adjunctions $f\dashv u$ and $v\dashv g$ compose to give an adjunction
$fv\dashv gu$. Similarly, the adjunctions $\overline{f}\dashv\overline{u}$ and
$\overline{v}\dashv\overline{g}$ compose to give an adjunction
$\overline{v}\overline{f}\dashv\overline{u}\,\overline{g}$. The result now
follows immediately from the lemma.  
\endproof

\begin{theorem}
Let \K be a 2-category with Eilenberg-Moore objects for monads and comonads
and splittings for idempotent 2-cells. A commutative square 
$$
\xymatrix{
P \ar[r]^{\overline{v}} \ar[d]_{\overline{u}} & L \ar[d]^u \\
M \ar[r]^-v & K}
$$
is in the image of the embedding $\WDL(\K)\to\K^{\two\x\two}$ if and only if
$u$ and $\overline{u}$ are monadic, $v$ and $\overline{v}$ are comonadic, and
the mate $\mone:fv\to\overline{v}\overline{f}$ (under the adjunctions $f\dashv
u$ and ${\overline f}\dashv {\overline u}$) has a section $\monesec$.
\end{theorem}

\proof
Since all these notions are preserved and jointly reflected by the
representable functors $\K(X,-):\K\to\Cat$, it will suffice to prove the
theorem in the case $\K=\Cat$. Suppose that $u$ and $\overline{u}$ are
monadic,  while $v$ and $\overline{v}$ are comonadic; this is certainly
necessary for the square to be in the image. Denote the various adjoints as
above, let $t$ be the monad induced by $f\dashv u$ and $c$ the comonad induced
by $v\dashv g$.  

Assume that the square is the image of some weak distributive law
$(K,t,c,\lambda)$ and take an object $(a,\gamma)$ of $M=K^c$. Then
$fv(a,\gamma)=(ta,\mu a)$, while $\overline{v}\overline{f}(a,\gamma)$ is
obtained by splitting (that is, taking the image of) the idempotent  
$$
\xymatrix{
ta \ar[r]^{t\gamma} & tca \ar[r]^{\lambda a} & cta \ar[r]^{\varepsilon ta} & ta.}
$$
By the last of the diagrams defining a weak distributive law, this idempotent
is in fact a morphism of $t$-algebras. The epimorphism appearing in the
splitting is just the component at $(a,\gamma)$ of
$\mone:fv\to\overline{v}\overline{f}$, and so we can take the other half of
the splitting as $\monesec$. This proves that any square in the image
satisfies the stated conditions.  

Suppose conversely that our square satisfies the stated condition. We have a
monad $t$ and a comonad $c$ on $K$. We construct $\lambda:tc\to ct$ as the
composite 
$$
\xymatrix{
ufvg \ar[r]^{u\mone g} & u\overline{v}\overline{f}g \ar@{=}[r] &
v\overline{u}\overline{f}g \ar[r]^{v\overline{u}\mthree} &
v\overline{u}\,\overline{g}f \ar[r]^{v\mtwo f} & 
vguf }
$$
where $\mthree:\overline{f}g\to\overline{g}f$ is the mate of $\monesec$ under
the adjunctions $v\dashv g$ and $\overline{v}\dashv\overline{g}$ (or
alternatively, the mate of  $\mtworetr$  under $f\dashv u$ and
$\overline{f}\dashv\overline{u}$). We must show that $\lambda$ is a weak
distributive law, and that the category of mixed $(K,t,c,\lambda)$-algebras is 
$P$. 

In the diagram 
$$
\xymatrix @R=10pt{
ufufvg \ar[r]^{ufu\mone g} \ar[ddddd]_{u\varepsilon fvg} & 
ufu\overline{v}\overline{f}g \ar@{=}[r]
\ar[ddddd]_{u\varepsilon\overline{v}\overline{f}g} & 
ufv\overline{u}\overline{f}g \ar[r]^{ufv\overline{u}\mthree}
\ar[dd]_{u\mone\overline{u}\overline{f}g} & 
ufv\overline{u}\,\overline{g}f \ar[r]^{ufv\mtwo f} \ar@{=}[d] & 
ufvguf \ar[d]^{u\mone guf} \ar[dl]^{ufv\mtworetr f} \\
&&&ufv\overline{u}\,\overline{g}f \ar[d]_{u\mone\overline{u}\,\overline{g}f} & 
u\overline{v}\overline{f}guf \ar@{=}[d] \\
&& u\overline{v}\overline{f}\overline{u}\overline{f}g \ar@{=}[d] 
\ar[dddl]_{u\overline{v}\,\overline{\varepsilon}\overline{f}g} & 
u\overline{v}\overline{f}\overline{u}\,\overline{g}f \ar@{=}[d] & 
v\overline{u}\overline{f}guf \ar[d]^{v\overline{u}\mthree uf}
\ar[dl]^{v\overline{u}\overline{f}\mtworetr f}  \\ 
&& v\overline{u}\overline{f}\overline{u}\overline{f}g
\ar[dd]_{v\overline{u}\,\overline{\varepsilon}\overline{f}g} & 
v\overline{u}\overline{f}\overline{u}\,\overline{g}f
\ar[dd]_{v\overline{u}\,\overline{\varepsilon}\,\overline{g}f} &
v\overline{u}\,\overline {g}fuf \ar[d]^{v\mtwo fuf}
\ar[ddl]_(0.6){v\overline{u}\,\overline{g}\varepsilon f} \\ 
&&&& vgufuf \ar[d]^{vgu\varepsilon f} \\
ufvg \ar[r]^-{u\mone g} & 
u\overline{v}\overline{f}g \ar@{=}[r] & 
v\overline{u}\overline{f}g \ar[r]^-{v\overline{u}\mthree} & 
v\overline{u}\,\overline{g}f \ar[r]^-{v \mtwo f} & vguf }$$
\noindent 
the triangle in the top right corner commutes since $\hat{\sigma}\sigma=1$,
the remaining regions commute by naturality or by mateship relations, and so
the exterior commutes, giving compatibility of $\lambda$ with
$\mu$. Compatibility with $\delta$ is similar. 

In the diagram
$$
\xymatrix @R=10pt{
vg \ar[r]^{v\eta g} \ar[d]_{\eta vg} \ar[dddr]^{v\overline{\eta}g} & 
vgvg \ar[r]^{vg\eta vg} \ar@/_1pc/[rrr]_{vgv\overline{\eta}g} & 
vgufvg \ar[r]^{vgu\pi g} & vgu\overline{v}\overline{f}g \ar@{=}[r] &
vgv\overline{u}\overline{f}g \ar[d]^{vgv\overline{u}\mthree} \\
ufvg \ar[dd]_{u\mone g} &&&& vgv\overline{u}\,\overline{g}f \ar[d]^{vgv\mtwo f} \\
&&&& vgvguf \ar[d]^{vg\varepsilon uf} \\
u\overline{v}\overline{f}g \ar@{=}[r]
& v\overline{u}\overline{f}g  \ar[r]^-{v\overline{u}\mthree}
\ar[uuurrr]^{v\eta\overline{u}\overline{f}g} & 
v\overline{u}\,\overline{g}f \ar[r]^-{v \mtwo f}& 
vguf  \ar@{=}[r] \ar[ur]^{v\eta guf} & vguf }
$$
the triangle in the bottom right corner commutes by one of the triangle
equations for the adjunction $v\dashv g$, the triangular region on the left
and the narrow region at the top commute by mateship relations, the remaining
regions commute by naturality, and so the exterior commutes, giving
compatibility of $\lambda$ with $\eta$. Compatibility with $\varepsilon$ is
similar. 

Finally we should observe that the monad induced by the adjunction
$\overline{f}\dashv\overline{u}$ is precisely the weakly lifted monad
$\overline{t}$, and so $P$ is the category of $\overline{t}$-algebras, which
in turn is the category of mixed $(K,t,c,\lambda)$-algebras. 
\endproof

\begin{remark}
Among the weak mixed distributive laws are the actual (non-weak) mixed
distributive laws. They correspond to the case where
$\pi:fv\to\overline{v}\overline{f}$ is not just a split epimorphism, but 
invertible. Since $fv\dashv gu$ and
$\overline{v}\overline{f}\dashv\overline{u}\,\overline{g}$, this is of course 
equivalent to $\sigma:\overline{u}\,\overline{g}\to gu$ being invertible.  
\end{remark}

\section{Consequences}

In this final section we draw a few simple consequences from what has gone
before. We continue to suppose that \K has Eilenberg-Moore objects for monads 
and comonads and splittings for idempotent 2-cells.

We shall use several times the following standard result:
\begin{lemma}\label{lem:adjoints}
For any commutative square 
$$
\xymatrix{
A \ar[r]^{G} \ar[d]_{H'} & B \ar[d]^{H} \\
C \ar[r]^-{G'} & D }
$$
with $H$ and $H'$ fully faithful, if $G'$ has a right adjoint $R'$ and $R'H$
lands in $A$, in the sense that $R'H\cong H'R$ for some $R$, then this $R$
gives a right adjoint to $G$. In particular, if $H$ is fully faithful, and
$HG$ has a right adjoint then so does $G$. 
\end{lemma}

\proof 
We have natural isomorphisms
$B(Ga,b)\cong D(HGa,Hb)=D(G'H'a,Hb)\cong C(H'a,R'Hb)\cong C(H'a,H'Rb)\cong
A(a,Rb)$.
\endproof

First of all, observe that the composite of the embedding 
$H:\WDL(\K)\to\K^{\two\x\two}$ and $Y:\K\to\WDL(\K)$ is the diagonal
2-functor $\Delta:\K\to\K^{\two\x\two}$ sending an object $X$ to the square
$$
\xymatrix{
X \ar@{=}[r] \ar@{=}[d] & X \ar@{=}[d] \\
X \ar@{=}[r] & X. }
$$
We constructed the embedding $H$ using the existence of Eilenberg-Moore
objects, but conversely, from the existence of a fully faithful
$H:\WDL(\K)\to\K^{\two\x\two}$ satisfying $HY=\Delta$, we may deduce the
existence of a right adjoint to $Y$ by Lemma \ref{lem:adjoints},
for $\Delta$ certainly has a right adjoint, sending a square to domain of the
diagonal (the top left corner in the way we have been drawing squares). 

\begin{proposition}
There is a fully faithful 2-functor $J:\Cmd(\K)\to\WDL(\K)$  sending a comonad
$(K,c)$ to $(K,1,c,1)$, and this 2-functor has a right 2-adjoint. 
\end{proposition}

\proof
Here we use the embedding $H_c:\Cmd(\K)\to\K^\two$ sending a comonad $(K,c)$
to the forgetful $v:K^c\to K$. We also use the fully faithful map
$J':\K^\two\to\K^{\two\x\two}$ sending a morphism $v:M\to K$ to the square  
$$
\xymatrix{
M \ar[r]^{v} \ar@{=}[d] & K \ar@{=}[d] \\
M \ar[r]^-v & K. }
$$
The composite $J'H_c$ clearly lands in the image of $H$, and so we obtain a
fully faithful $J$ with $HJ=J'H_c$. Now  $J'$  has a right adjoint sending a
square 
$$
\xymatrix{
P \ar[r]^{\overline{v}} \ar[d]_{\overline{u}} & L \ar[d]^u \\
M \ar[r]^-{v} & K}
$$
to $\overline{v}:P\to L$, and this map is comonadic if the original square is
in the image of $\WDL(\K)$, thus we obtain by Lemma \ref{lem:adjoints}
the desired right adjoint. 
\endproof

Dually we have 

\begin{proposition}
There is a fully faithful 2-functor $J_*:\Mnd(\K)\to\WDL(\K)$  sending a 
monad $(K,t)$ to $(K,t,1,1)$, and this 2-functor has a right
2-adjoint. \endproof  
\end{proposition}
 
Finally we observe that the diagram 
$$
\xymatrix{
\K \ar[r]^{I} \ar[d]_{I_*} & \Mnd(\K) \ar[d]^{J_*} \\
\Cmd(\K) \ar[r]^-{J} & \WDL(\K) }
$$
of fully faithful 2-functors commutes, with the diagonal being what we have
called $Y$. Thus the corresponding diagram of right adjoints commutes up to
natural isomorphism.

\end{document}